\documentclass[11pt,a4paper]{article}
\usepackage{a4wide,amsfonts,amsmath,latexsym,amssymb,euscript,graphicx,units,mathrsfs}

\usepackage{amsthm, amsfonts, amssymb,eufrak}

\usepackage{float}
\newfloat{figure}{H}{lof}
\floatname{figure}{\figurename}

\DeclareMathAlphabet{\eufrak}{U}{}{}{}  
\SetMathAlphabet\eufrak{normal}{U}{euf}{m}{n}
\SetMathAlphabet\eufrak{bold}{U}{euf}{b}{n}

\usepackage{amsmath, amsthm, amsfonts, amssymb}

\numberwithin{equation}{section}

\def\real{{\mathord{{\rm I\kern-2.8pt R}}}}        
\def\inte{{\mathord{{\rm I\kern-2.8pt N}}}}
\def\PP{{\mathord{{\rm I\kern-2.8pt P}}}}

\def\real{{\mathord{\mathbb R}}}

\def\inte{{\mathord{\mathbb N}}}

\def\R{\right}
\def\L{\left}

\def\P{\mathbb{P}}
\def\E{\mathop{\hbox{\rm I\kern-0.20em E}}\nolimits}

\newtheorem{prop}{Proposition}[section]

\newtheorem{theorem}[prop]{Theorem}

\allowdisplaybreaks

\begin{document}

\vspace*{-2.1cm} 

\hfill{Probabilit\'es}
\bigskip 

\begin{center}
{\bf 
\Large
On the orthogonal component of BSDEs in a Markovian setting
} 
\normalsize
\end{center}

\begin{center}
\large{\footnote{Humboldt-Universit\"at zu Berlin\\
Institut f\"ur Mathematik\\
Unter den Linden 6\\
10099 Berlin Germany\\ {\tt areveill@mathematik.hu-berlin.de}}{Anthony R\'eveillac}}
\normalsize
\end{center}
{\bf Abstract - }
{ 
\small 
In this Note we consider a quadratic backward stochastic differential equation (BSDE) driven by a continuous martingale $M$ and whose generator is a deterministic function. We prove (in Theorem \ref{theorem:main}) that if $M$ is a strong homogeneous Markov process and if the BSDE has the form \eqref{BSDE} then the unique solution $(Y,Z,N)$ of the BSDE is reduced to $(Y,Z)$, \textit{i.e.} the orthogonal martingale $N$ is equal to zero showing that in a Markovian setting the "usual" solution $(Y,Z)$ has not to be completed by a strongly orthogonal even if $M$ does not enjoy the martingale representation property.
\normalsize 
} 

\begin{center}
{\bf Sur la composante orthogonale d'une EDSR dans un contexte markovien 
} 
\normalsize
\end{center}
{\em {\bf R\'esum\'e - }} 
{
\small 
Dans cette Note nous consid\'erons une \'equation diff\'erentielle stochastique r\'etrograde (EDSR) de g\'en\'erateur d\'eterministe et quadratique dirig\'ee par une martingale continue $M$. Nous prouvons (dans le Th\'eor\`eme \ref{theorem:main}) que si $M$ est un processus de Markov homog\`ene fort et si l'EDSR est de la forme \eqref{BSDE} l'unique solution $(Y,Z,N)$ de l'EDSR se r\'eduit \`a $(Y,Z)$, \textit{i.e.} la martingale orthogonale $N$ vaut z\'ero. Cela prouve que dans un contexte markovien la solution "habituelle" $(Y,Z)$ n'a pas \`a \^etre compl\'et\'ee par une martingale fortement orthogonale m\^eme si $M$ ne poss\`ede pas la propri\'et\'e de repr\'esentation martingale. 
\normalsize 	
} 

\bigskip
\baselineskip0.5cm 
\noindent {\bf Version fran\c caise abr\'eg\'ee} 
\\
Dans cette Note nous consid\'erons une \'equation diff\'erentielle stochastique r\'etrograde (EDSR) dirig\'ee par une martingale continue $M$, de g\'en\'erateur quadratique $f$ et admettant $F(X_T)$ pour condition terminale o\`u $F:\real\to\real$ d\'enote une fonction d\'eterministe suffisamment r\'eguli\`ere et $X$ l'unique solution forte d'une \'equation diff\'erentielle stochastique (EDS) \'egale\-ment dirig\'ee par $M$. Dans ce contexte il a \'et\'e d\'emontr\'e dans \cite{ElKarouiHuang} and \cite{Morlais} qu'il existe un unique triplet $(Y,Z,N)$ solution de l'EDSR consid\'er\'ee o\`u $Y$ est un processus stochastique born\'e, $Z$ un processus pr\'evisible de carr\'e int\'egrable et $N$ une martingale fortement orthogonale \`a $M$. Puisque nous ne supposons pas que $M$ poss\`ede la propri\'et\'e de repr\'esentation martingale, la solution habituelle $(Y,Z)$ doit \textit{a priori} \^etre compl\'et\'ee par une martingale $N$ fortement orthogonale \`a $M$. Si le g\'en\'erateur $f$ est suppos\'e Lipschitz, les auteurs de \cite{ElKarouiHuang} obtiennent la solution de l'EDSR \eqref{BSDE} \textit{via} une it\'eration de Picard de la forme \eqref{eq.Picard}. Notons que la troisi\`eme compostante de la solution, la martingale orthogonale $N$ est "statique" lors de cette it\'eration.\\\\
\noindent 
L'objet de cette Note est de d\'emontrer que dans un contexte markovien (\textit{i.e.} avec une condition terminale comme expos\'ee plus haut et un g\'en\'erateur d\'eterministe d\'ependant uniquement de $y$ et $z$) la solution $(Y,Z,N)$ se r\'eduit au couple $(Y,Z)$ autrement dit, la composante orthogonale $N$ est nulle m\^eme si la propri\'et\'e de repr\'esentation martingale n'est pas v\'erifi\'ee pour $M$. Afin de simplifier la preuve du r\'esultat principal (Th\'eor\`eme \ref{theorem:main}) nous consid\'erons une diffusion $X$ de d\'erive nulle et toutes les \'equations mises en jeu sont uni-dimensionnelles (le cas d'un g\'en\'erateur d\'ependant de $(X,M)$ fera l'objet d'un travail futur). Ce r\'esultat permettera (dans un travail en pr\'eparation) de simplifier l'\'etude des propri\'et\'es des EDSR quadratiques de la forme \eqref{BSDE} comme en particulier donner une preuve de diff\'erentiabilit\'e par rapport aux param\`etres initiaux $(x,m)$ (voir \eqref{SDE}) sans l'hypoth\`ese additionnelle (MRP) (\textit{c.f.} \cite[Section 4.2]{ImkellerReveillacRichter}) utilis\'ee dans \cite[Theorem 4.6]{ImkellerReveillacRichter}. 
\vspace{-0.03cm} 
\noindent 

\section{Preliminaries}
\label{section1}

Let $M:=(M_t)_{t\in [0,T]}$ be a real-valued continuous square integrable martingale with respect to a continuous filtration $(\mathcal{F}_t)_{t\in [0,T]}$ both defined on a probability space $(\Omega,\mathcal{F},\P)$. Assume that $M$ is an homogeneous strong Markov process with respect to $(\mathcal{F}_t)_{t\in [0,T]}$. For $(t,m)$ in $[0,T]\times \real$ we denote by $M^{t,m}$ the process defined as $ M_s^{t,m}:=m+M_s-M_t, \quad s \in [t,T].$ Let $C:=(C_t)_{t\in[0,T]}$ be the $(\mathcal{F}_t)_{t\in [0,T]}$-predictable, increasing process defined by $ C_t:=\arctan(\langle M, M \rangle_t), \; t \in [0,T].$
On this filtered probability space we also consider a stochastic process $X^{t,x,m}:=(X_s^{t,x,m})_{s\in [t,T]}$ defined as the unique strong solution of the following one-dimensional stochastic differential equation
\begin{equation}
\label{SDE}
X_s^{t,x,m}=x+\int_t^s \sigma(X_r^{t,x,m},M_r^{t,m}) dM_r, \quad s \in[t,T], \; t \in [0,T]
\end{equation}
where $\sigma:\real\times\real \to \real$ is deterministic, of class $C^2(\real\times \real)$ with locally Lipschitz partial derivatives and such that there exists a positive constant $k$ satisfying $ \vert \sigma(x_1,m_1)-\sigma(x_2,m_2) \vert \leq k \vert x_1-x_2 \vert, \quad \forall (x_1,x_2,m_1,m_2) \in \real^4.$
Let us finally introduce the object of interest of this Note that is the following backward stochastic differential equation (BSDE) coupled with the forward process $X^{t,x,m}$ as
\begin{align} 
\label{BSDE}
Y_s^{t,x,m}=&F(X_T^{t,x,m})-\int_t^T Z_r^{t,x,m} dM_r + \int_t^T f(r,Y_r^{t,x,m},Z_r^{t,x,m}) dC_r-\int_t^T dN_r^{t,x,m}\nonumber\\
&+ \frac{\kappa}{2} \int_t^T d\langle N^{t,x,m},N^{t,x,m}\rangle_r.
\end{align}
where $F:\real\to\real$ is a bounded deterministic function of class $\mathcal{C}^2(\real)$ with bounded derivatives. The generator $f: [0,T] \times \real \times \real \to \real$ is assumed to be $\mathcal{B}([0,T]) \otimes \mathcal{B}(\real) \otimes \mathcal{B}(\real)$-measurable where $\mathcal{B}(\real)$ is for the Borel $\sigma$-filed on $\real$ (so that $f(r,x,m)$ is deterministic for non-random $(r,x,m)$ in $[0,T]\times\real^2$) and is such that there exists a deterministic constant $c$ satisfying $ \sup_{r \in [0,T]} \vert f(r,0,0) \vert \leq c.$ 
We assume in addition that the generator $f$ is quadratic in $z$ and Lipschitz in $y$. The typical example being  when $f$ is of the form $f(s,y,z)=l(s,y)+\eta \vert z \vert^2$ where $\eta$ is a fixed constant and $l$ is Lipschitz in $y$ (the more general "quadratic" assumptions can be found for example in \cite{ImkellerReveillacRichter}). We recall that in this setting, it is shown in \cite{Morlais} that there exists a unique triple $(Y^{t,x,m}, Z^{t,x,m}, N^{t,x,m}) \in \mathcal{S}^\infty\times L^2(d \langle M, M \rangle \otimes d\P)\times \mathcal{M}^2$ where $\mathcal{S}^\infty$ is the space of bounded and continuous $(\mathcal{F}_t)_t$-adapted processes, $L^2(d \langle M, M \rangle \otimes d\P)$ denotes the space of square integrable $(\mathcal{F}_t)_t$-predictable processes and $\mathcal{M}^2$ the space of square integrable $(\mathcal{F}_t)_t$-martingales $N$ strongly orthogonal to $M$ (\textit{i.e.} $\langle M, N \rangle=0$). We also mention that these processes are real-valued.
We finally stress that all the conditions and assumptions previously mentioned will be assumed to hold in the rest of this Note and that $K$ denotes a constant which can differ from one line to another. We conclude this section by recalling some important facts. First let us mention that only the couple $(X,M)$ is an homogeneous strong Markov process.

\begin{theorem}(\cite[Theorem (8.11)]{CinlarJacodProtterSharpe} or \cite[V. Theorem 35]{Protter})
The process $(X_s^{t,x,m},M_s^{t,m})_{s\in [t,T]}$ is an homogeneous strong Markov process for the filtration $(\mathcal{F}_t)_{t\in [0,T]}$. If in addition $M$ is assumed to enjoy the independent increments property then the stochastic process $(X_s^{t,x,m})_{s\in [t,T]}$ is a strong Markov process.
\end{theorem} 

\noindent
The Markov property of the couple $(X,M)$ transfers to the solution of \eqref{BSDE} and \eqref{eq.first}.

\begin{theorem}(\cite[Proposition 3.2, Theorem 3.4]{ImkellerReveillacRichter})
\label{prop:Markovprop}
There exist two deterministic functions $u,v:[0,T] \times \real^2 \to \real$, $\mathcal{B}([0,T]) \otimes \mathcal{B}(\real^2)$ such that $(Y^{t,x,m},Z^{t,x,m})$ in \eqref{BSDE} and \eqref{eq.first} satisfy
$$ Y_s^{t,x,m}=u(s,X_s^{t,x,m},M_s^{t,m}), \quad Z_s^{t,x,m}=v(s,X_s^{t,x,m},M_s^{t,m}) \sigma(s,X_s^{t,x,m},M_s^{t,m}), \quad s \in [t,T]$$
where $\mathcal{B}_e(\real^2)$ is the $\sigma$-field on $\real^2$ generated by functions $(x,m) \mapsto \E\L[\phi(s,X_s^{t,x,m},M_s^{t,m}) dC_s\R]$ with $\phi:\Omega \times [0,T] \times \real^2 \to \real$ a continuous bounded function.
\end{theorem}
\noindent
Finally we will use the following \footnote{Note that this result (\cite[Theorem 4.6]{ImkellerReveillacRichter}) has been proved under an additional technical assumption (MRP) with $f$ a quadratic generator. Since the generator in equation \eqref{eq.first} is very simple, using only an exponential change we can apply the computations realized in \cite[Theorem 4.6]{ImkellerReveillacRichter} without assuming the hypothesis (MRP). The full proof of this fact will be presented in a paper in preparation.}{property} for the solution of the BSDE \eqref{eq.first}.

\begin{theorem}(Particular case of \cite[Theorem 4.6]{ImkellerReveillacRichter})
\label{prop:Diff} 
The map $ (x,m)\mapsto Y^{1,t,x,m} $ is of class $\mathcal{C}^1(\real^2)$ $\P$-a.s. where $Y^{1,t,x,m}$ is as in \eqref{eq.first} below.
\end{theorem}

\section{Main result}

We are now ready to state and prove the main result of this Note.

\begin{theorem}
\label{theorem:main}
Assume that assumptions of Section \ref{section1} are in force then $N^{t,x,m}$ in \eqref{BSDE} is equal to zero and equation \eqref{BSDE} becomes
$$ Y_s^{t,x,m}=F(X_T^{t,x,m})-\int_t^T Z_r^{t,x,m} dM_r + \int_t^T f(r,Y_r^{t,x,m},Z_r^{t,x,m}) dC_r.$$
\end{theorem}

\begin{proof}
First note that it is enough to assume that the generator $f$ is Lipschitz in $(y,z)$. Indeed, in \cite[Theorems 2.5-2.6]{Morlais}, the existence and uniqueness of the solution $(Y^{t,x,m},Z^{t,x,m},N^{t,m})$ of the BSDE \eqref{BSDE} is given as a limit of solutions of Lipschitz BSDEs. As a consequence, $N^{t,x,m}$ is equal to zero in \eqref{BSDE} if the orthogonal martingale parts $N$ in the approximating Lipschitz BSDEs vanish. So assume $f$ to be Lipschitz in $(y,z)$. In \cite{ElKarouiHuang} the authors show that the unique solution of \eqref{BSDE} is obtained as the limit of the following Picard iteration:
\begin{align}
\label{eq.Picard}
Y_s^{0,t,x,m}&=Z_s^{0,t,x,m}=0,\nonumber\\
Y_s^{k+1,t,x,m}&=F(X_T^{t,x,m})- \int_s^T Z_r^{k,t,x,m} dM_r +\int_s^T f(r,Y_r^{k,t,x,m},Z_r^{k,t,x,m}) dC_r\nonumber\\
&- \int_s^T dN_r^{t,x,m} + \kappa \int_s^T d\langle N^{t,x,m}, N^{t,x,m}\rangle_r, \quad k\geq 0. 
\end{align} 
Note that $N^{t,x,m}$ is not part of the iteration (we refer to \cite[Proof of Theorem 6.1]{ElKarouiHuang} for more details). This remark leads to the main idea of the proof. Indeed, consider the first iteration, \textit{i.e.} $(Y^{1,t,x,m}, Z^{1,t,x,m}, N^{t,x,m})$ defined by 
\begin{equation}
\label{eq.first}
Y_s^{1,t,x,m}=F(X_T^{t,x,m})- \int_s^T Z_r^{1,t,x,m} dM_r +\int_s^T f(r,0,0) dC_r- \int_s^T dN_r^{t,x,m} + \kappa \int_s^T d\langle N^{t,x,m}, N^{t,x,m}\rangle_r.
\end{equation}
By the a priori estimates obtained in \cite[Proposition 6.3]{ElKarouiHuang} the triplet $(Y^{1,t,x,m}, Z^{1,t,x,m}, N^{t,x,m})$ is unique. As a consequence if we show that $N^{t,x,m}=0$ in equation \eqref{eq.first} then the Theorem is proved since $(Y^{k,t,x,m}, Z^{k,t,x,m}, N^{t,x,m})$ converges to the unique solution of \eqref{BSDE}. The rest of the proof is devoted to this fact.\\\\
\noindent 
Since $Y^{1,t,x,m}$ is $(\mathcal{F}_\cdot)$-adapted it holds by Markov property that
$$ Y_s^{1,t,x,m}=g(s,X_s^{t,x,m},M_s^{t,m}), \; \textrm{ with } g(s,x,m):=\E\left[F(X_{T-s}^{t,x,m}) -\int_s^T f(r,0,0) dC_r\right].$$  
In addition, Proposition \ref{prop:Diff} applied to \eqref{eq.first} gives that the application $(x,m) \mapsto g(t,x,m)$ is of class $\mathcal{C}^1(\real\times\real)$ for every $t$. 
We mimic a technique given in \cite{ImkellerReveillacRichter} and compute $\langle Y^{1,t,x,m}, N^{1,t,x,m} \rangle_s$ for $s\geq t$. Let $\pi^{(n)}:=\{t=t_0^{(n)} \leq t_1^{(n)} \leq \cdots \leq t_N^{(n)}=s\}$ be a family of subdivisions of $[t,s]$ whose mesh $\vert \pi^{(n)} \vert$ tends to zero as $n$ goes to the infinity. For sake of simplicity the superscript $(n)$ will be omitted in the following computations.  
\begin{eqnarray}
\label{bracket1}
\langle Y^{1,t,x,m}, N^{t,x,m} \rangle_s&=&\langle g(\cdot,X_\cdot^{t,x,m},M_\cdot^{t,x,m}), N^{t,x,m} \rangle_s\nonumber\\ 
&\overset{\P}{=}&\lim_{n \to \infty} \sum_{j=1}^r (g(t_{j+1},X_{t_{j+1}}^{t,x,m},M_{t_{j+1}}^{t,x,m})-g(t_j,X_{t_j}^{t,x,m},M_{t_j}^{t,x,m})) \Delta_j N^{t,x,m}\nonumber\\
&\overset{\P}{=}& \lim_{n \to \infty} \sum_{j=1}^r\bigg[ (g(t_j,X_{t_{j+1}}^{t,x,m},M_{t_{j+1}}^{t,x,m})-g(t_j,X_{t_j}^{t,x,m},M_{t_j}^{t,x,m})) \Delta_j N^{t,x,m}\nonumber\\
&&+(g(t_{j+1},X_{t_{j+1}}^{t,x,m},M_{t_{j+1}}^{t,x,m})-g(t_j,X_{t_{j+1}}^{t,x,m},M_{t_{j+1}}^{t,x,m})) \Delta_j N^{t,x,m} \bigg].
\end{eqnarray}
We consider the two sumands above separately. For the first part we follow a technique used in \cite{ImkellerReveillacRichter} and apply the mean theorem. Let $\bar{M}_j$ (respectively $\bar{X}_j$) below a random point between $M_{t_j}^{t,x,m}$ and $M_{t_{j+1}}^{t,x,m}$ (resp. $X_{t_j}^{t,x,m}$ and $X_{t_{j+1}}^{t,x,m}$) in the computations below. We have  
\begin{eqnarray}
\label{bracket2}
&&\sum_{j=1}^r (g(t_j,X_{t_{j+1}}^{t,x,m},M_{t_{j+1}}^{t,x,m})-g(t_j,X_{t_j}^{t,x,m},M_{t_j}^{t,x,m})) \Delta_j N^{t,x,m}\nonumber\\
&=&\sum_{j=1}^r (g(t_j,X_{t_{j+1}}^{t,x,m},M_{t_{j+1}}^{t,x,m})-g(t_j,X_{t_j}^{t,x,m},M_{t_{j+1}}^{t,x,m})) \Delta_j N^{t,x,m}\nonumber\\
&&+\sum_{j=1}^r (g(t_j,X_{t_j}^{t,x,m},M_{t_{j+1}}^{t,x,m})-g(t_j,X_{t_j}^{t,x,m},M_{t_j}^{t,x,m})) \Delta_j N^{t,x,m}\nonumber\\
&=&\sum_{j=1}^r \bigg[\partial_2 g(t_j,X_{t_j}^{t,x,m},M_{t_j}^{t,x,m}) \Delta_j X \Delta_j N^{t,x,m}+ \partial_3 g(t_j,X_{t_j}^{t,x,m},M_{t_j}^{t,x,m}) \Delta_j M \Delta_j N^{t,x,m} +R_{j,r}\bigg] \quad\quad\quad{}
\end{eqnarray}
where $R_{j,r}$ is defined as
\begin{eqnarray*}
R_{j,r}&:=& (\partial_2 g(t_j,\bar{X}_j,M_{t_{j+1}}^{t,x,m}-\partial_2 g(t_j,X_{t_j}^{t,x,m},M_{t_j}^{t,x,m})) \Delta_j X \Delta_j N^{t,x,m}\\
&&+ (\partial_3 g(t_j,X_{t_j}^{t,x,m},\bar{M}_j) -\partial_3 g(t_j,X_{t_j}^{t,x,m},M_{t_j}^{t,x,m})) \Delta_j M \Delta_j N^{t,x,m}.
\end{eqnarray*}
Since $(x,m) \mapsto g(s,x,m)$ is of class $\mathcal{C}^1$ for every $s$ in $[0,T]$ the remainder term $\sum_{j=0}^r R_{j,r}$ as $r$ goes to infinity (we refer to \cite[Proof of (5.13)]{ImkellerReveillacRichter} for the complete justifications). Then it follows using \eqref{bracket2} that
\begin{eqnarray*}
&&\lim_{r\to\infty} \sum_{j=1}^r (g(t_j,X_{t_{j+1}}^{t,x,m},M_{t_{j+1}}^{t,x,m})-g(t_j,X_{t_j}^{t,x,m},M_{t_j}^{t,x,m})) \Delta_j N^{t,x,m}\\
&=&\langle \int_t^\cdot \partial_2 g(r,X_r^{t,x,m},M_r^{t,x,m}) \sigma(r,X_r^{t,x,m},M_r^{t,x,m})+ \partial_3 g(r,X_r^{t,x,m},M_r^{t,x,m}) dM_r, N_\cdot^{t,x,m} \rangle_s=0
\end{eqnarray*}
by strong orthogonality between $M$ and $N$.
As a consequence, relation \eqref{bracket1} reduces to
\begin{equation}
\label{bracket5}
\langle Y^{1,t,x,m}, N^{t,x,m} \rangle_s\overset{\P}{=} \lim_{n \to \infty} \sum_{j=1}^r (g(t_{j+1},X_{t_{j+1}}^{t,x,m},M_{t_{j+1}}^{t,x,m})-g(t_j,X_{t_{j+1}}^{t,x,m},M_{t_{j+1}}^{t,x,m})) \Delta_j N^{t,x,m}.
\end{equation}
We have that
\begin{eqnarray*}
&&\L\vert \sum_{j=1}^r (g(t_{j+1},X_{t_{j+1}}^{t,x,m},M_{t_{j+1}}^{t,x,m})-g(t_j,X_{t_{j+1}}^{t,x,m},M_{t_{j+1}}^{t,x,m})) \Delta_j N^{t,x,m} \R\vert^2\\
&\leq&\sum_{j=1}^n \left\vert g(t_{j+1},X_{t_{j+1}}^{t,x,m},M_{t_{j+1}}^{t,x,m})-g(t_j,X_{t_{j+1}}^{t,x,m},M_{t_{j+1}}^{t,x,m})\right\vert^2 \times \sum_{j=1}^n \vert \Delta_j N \vert^2\\
&=&\sum_{j=1}^n \left\vert\E\left[F(X_{T-t_{j+1}}^{0,X_{t_{j+1}}^{t,x,m},M_{t_{j+1}}^{t,m}})-F(X_{T-t_j}^{0,X_{t_{j+1}}^{t,x,m},M_{t_{j+1}}^{t,m}})-\int_{t_j}^{t_{j+1}} f(r,0,0) dC_r\right]\right\vert^2 \times \sum_{j=1}^n \vert \Delta_j N \vert^2\\
&\leq& 2 \left[\sum_{j=1}^n \left\vert\E\left[F(\tilde{X}_{T-t_{j+1}})-F(\tilde{X}_{T-t_j}) \right]\right\vert^2 + \sum_{j=1}^n \left\vert\E\left[\int_{t_j}^{t_{j+1}} f(r,0,0) dC_r\right]\right\vert^2 \right]\times \sum_{j=1}^n \vert \Delta_j N \vert^2\\
\end{eqnarray*}
where for simplicity of notations we set $\tilde{X}_s:=X_s^{0,X_{t_{j+1}}^{t,x,m},M_{t_{j+1}}^{t,m}}$.
Let $\bar{X}_j$ be a random point between $\tilde{X}_{T-t_{j+1}}$ and $\tilde{X}_{T-t_j}$. Writing $\E\L[F(\tilde{X}_{T-t_{j+1}})-F(\tilde{X}_{T-t_j})\R]$ as
\begin{eqnarray*}
&&\E\L[F(\tilde{X}_{T-t_{j+1}})-F(\tilde{X}_{T-t_j})\R]\\
&=& \E\L[F'(\tilde{X}_{T-t_{j+1}}) \L(\tilde{X}_{T-t_{j+1}}-\tilde{X}_{T-t_j}\R)\R] + \frac12 \E\L[F''(\bar{X}_j) \L\vert \tilde{X}_{T-t_{j+1}}-\tilde{X}_{T-t_j}\R\vert^2\R]\\
&=& \E\L[F'(\tilde{X}_{T-t_{j+1}}) \E\L[\tilde{X}_{T-t_{j+1}}-\tilde{X}_{T-t_j}\vert \mathcal{F}_{T-t_{j+1}}\R]\R] + \frac12 \E\L[F''(\bar{X}_j) \L\vert \tilde{X}_{T-t_{j+1}}-\tilde{X}_{T-t_j}\R\vert^2\R]\\
&\leq& K \E\L[\L\vert \tilde{X}_{T-t_{j+1}}-\tilde{X}_{T-t_j}\R\vert^2\R]
\end{eqnarray*}
and $f(r,0,0)$ as $ f(r,0,0)= \max\{f(r,0,0),0\}-\max\{-f(r,0,0),0\}$
it follows that 
\begin{eqnarray*}
&&\left\vert\sum_{j=1}^n \left(g(t_{j+1},M_{t_{j+1}})-g(t_j,M_{t_{j+1}})\right) \Delta_j N \right\vert^2\\
&\leq& K \left[\sum_{j=1}^n \left\vert E\left[ \vert \tilde{X}_{T-t_{j+1}}-\tilde{X}_{T-t_j}\vert^2 \right] \right\vert^2+ \sup_{r\in[0,T]} \vert f(r,0,0)) \vert \sum_{j=1}^n \left\vert\E\left[C_{t_{j+1}}-C_{t_j}\right]\right\vert^2 \right] \times \sum_{j=1}^n \vert \Delta_j N \vert^2\\
&\leq& K \left[\sum_{j=1}^n E\left[ \vert \tilde{X}_{T-t_{j+1}}-\tilde{X}_{T-t_j}\vert^4 \right] + \sup_{r\in[0,T]} \vert f(r,0,0)) \vert \sum_{j=1}^n \E\left[\vert M_{t_{j+1}}-M_{t_j} \vert^4\right] \right] \times \sum_{j=1}^n \vert \Delta_j N \vert^2\\
&\leq& K \left(E\left[\sum_{j=1}^\infty \vert \tilde{X}_{T-t_{j+1}}-\tilde{X}_{T-t_j}\vert^4 \right] + \E\left[ \sum_{j=1}^\infty \vert M_{t_{j+1}}-M_{t_j} \vert^4\right] \right) \times \sum_{j=1}^\infty \vert \Delta_j N \vert^2\\
&=&0
\end{eqnarray*} 
since the quartic variations of a martingale are zero. The previous computation and the equality \eqref{bracket5} entail that
\begin{equation}
\label{bracket6}
\langle Y^{1,t,x,m}, N^{t,x,m} \rangle_s\overset{\P}{=} 0. 
\end{equation}
On the other hand, the covariation $\langle Y^{1,t,x,m}, N^{t,x,m} \rangle_s$ in the BSDE \eqref{eq.first} equals to 
\begin{equation}
\label{bracket7}
\langle Y^{1,t,x,m}, N^{t,x,m} \rangle_s\overset{\P}{=} \langle N^{t,x,m}, N^{t,x,m} \rangle_s. 
\end{equation}
Hence relations \eqref{bracket6} and \eqref{bracket7} give that $N_s^{t,x,m}=N_0^{t,x,m}$ for every $s$ in $[t,T]$.
\end{proof}


%

\end{document}